
\magnification=1200
\input amstex.tex
\documentstyle{amsppt}

\hsize=12.5cm
\vsize=18cm
\hoffset=1cm
\voffset=2cm

\footline={\hss{\vbox to 2cm{\vfil\hbox{\rm\folio}}}\hss}
\nopagenumbers

\def\txt#1{{\textstyle{#1}}}
\baselineskip=13pt
\def\hf{{\textstyle{1\over2}}}
\def\a{\alpha}\def\b{\beta}
\def\d{{\,\roman d}}
\def\e{\varepsilon}
\def\f{\varphi}
\def\G{\Gamma}
\def\k{\kappa}
\def\s{\sigma}

\def\={\;=\;}

\def\zt{\zeta(\hf+it)}

\def\D{\Delta}
  
\def\R{\Re{\roman e}\,} \def\I{\Im{\roman m}\,} \def\s{\sigma}
\def\z{\zeta}

\def\ze{\zeta({1\over2} + it)}
\def\H{H_j^3({\txt{1\over2}})}  \def\={\,=\,}
\def\hf{{\textstyle{1\over2}}}
\def\txt#1{{\textstyle{#1}}}
\def\f{\varphi}
\def\z#1{|\zeta (\hf +i#1)|^4}

\font\tenmsb=msbm10
\font\sevenmsb=msbm7
\font\fivemsb=msbm5
\newfam\msbfam
\textfont\msbfam=\tenmsb
\scriptfont\msbfam=\sevenmsb
\scriptscriptfont\msbfam=\fivemsb
\def\Bbb#1{{\fam\msbfam #1}}

\def \NN {\Bbb N}
\def \CC {\Bbb C}
\def \RR {\Bbb R}
\def \ZZ {\Bbb Z}

\font\ff=cmr8
\def\txt#1{{\textstyle{#1}}}
\baselineskip=13pt
\def\hf{{\textstyle{1\over2}}}
\def\a{\alpha}\def\b{\beta}
\def\d{{\,\roman d}}
\def\e{\varepsilon}
\def\f{\varphi}
\def\G{\Gamma}
\def\k{\kappa}
\def\s{\sigma}

\def\={\;=\;}

\def\zt{\zeta(\hf+it)}

\def\D{\Delta}
  
\def\R{\Re{\roman e}\,} \def\I{\Im{\roman m}\,} \def\s{\sigma}
\def\z{\zeta}

\def\ze{\zeta({1\over2} + it)}
\def\H{H_j^3({\txt{1\over2}})}  \def\={\,=\,}
\font\teneufm=eufm10
\font\seveneufm=eufm7
\font\fiveeufm=eufm5
\newfam\eufmfam
\textfont\eufmfam=\teneufm
\scriptfont\eufmfam=\seveneufm
\scriptscriptfont\eufmfam=\fiveeufm
\def\mathfrak#1{{\fam\eufmfam\relax#1}}

\font\tenmsb=msbm10
\font\sevenmsb=msbm7
\font\fivemsb=msbm5
\newfam\msbfam
     \textfont\msbfam=\tenmsb
      \scriptfont\msbfam=\sevenmsb
      \scriptscriptfont\msbfam=\fivemsb
\def\Bbb#1{{\fam\msbfam #1}}

\def \NN {\Bbb N}
\def \CC {\Bbb C}

\def \RR {\Bbb R}
\def \ZZ {\Bbb Z}

  \def\rightheadline{{\hfil{\ff
On exponential sums with Hecke series at central points} \hfil\tenrm\folio}}

  \def\leftheadline{{\tenrm\folio\hfil{\ff
  Aleksandar Ivi\'c }\hfil}}
  \def\emptyheadline{\hfil}
  \headline{\ifnum\pageno=1 \emptyheadline\else
  \ifodd\pageno \rightheadline \else \leftheadline\fi\fi}

  \topmatter
\title
ON EXPONENTIAL SUMS WITH HECKE SERIES AT CENTRAL POINTS
\endtitle
\author   Aleksandar Ivi\'c \endauthor
\address{
Aleksandar Ivi\'c, Katedra Matematike RGF-a
Universiteta u Beogradu, \DJ u\v sina 7, 11000 Beograd,
Serbia (Yugoslavia).}
\endaddress
\email {\tt aivic\@matf.bg.ac.yu,
ivic\@rgf.bg.ac.yu} \endemail
\abstract
Upper bound  estimates  for the exponential sum
$$
\sum_{K<\k_j\le K'<2K} \a_j\H \cos\left(\k_j\log\left(
{4{\roman e}T\over \k_j}\right)\right) \qquad(T^\e \le K  \le T^{1/2-\e})
$$
are considered, where
$\a_j = |\rho_j(1)|^2(\cosh\pi\kappa_j)^{-1}$, and
$\rho_j(1)$ is the first Fourier coefficient of the Maass wave form
corresponding to the eigenvalue $\lambda_j = \k_j^2 + {1\over4}$
to which the Hecke series $H_j(s)$ is attached. The problem is transformed
to the estimation of a classical  exponential sum involving the binary
additive divisor problem. The analogous exponential sums with
$H_j(\hf)$ or $H_j^2(\hf)$ replacing $\H$ are also considered.
The above sum is  conjectured to be $\ll_\e K^{3/2+\e}$, which is
proved to be true in the mean square sense.
\endabstract
\endtopmatter
\title par Aleksandar Ivi\'c\endtitle

\def\DJ{\leavevmode\setbox0=\hbox{D}\kern0pt\rlap
{\kern.04em\raise.188\ht0\hbox{-}}D}

\head 1. Introduction
\endhead

\bigskip
The main purpose of this paper is to transform and estimate  exponential
sums of Hecke series at central points, namely  the sums
$$
S(K) = S(K;K',T) := \sum_{K<\k_j\le K'}\a_j\H\cos\left(\k_j\log{4{\roman e}T
\over\k_j}\right),\leqno(1.1)
$$
under the condition
$$
T^\e \le K < K' \le 2K \le T^{1/2-\e}.\leqno(1.2)
$$
Here and later $\e > 0$ denotes arbitrarily small constants,
not necessarily the same ones at each occurrence. The quantities
$\a_j, H_j(\hf)$ and $\k_j$ are connected with
the spectral theory of the non-Euclidean Laplacian.
For a comprehensive account of
spectral theory the reader is referred to Y. Motohashi's monograph [23],
and here we only briefly explain some basic notions.

\medskip
Let $\,\{\lambda_j = \kappa_j^2 + {1\over4}\}_{j=1}^\infty
 \,\cup\, \{0\}\,$ be the
eigenvalues (discrete spectrum) of the hyperbolic Laplacian
$$
\Delta=-y^2\left({\left({\partial\over\partial x}\right)}^2 +
{\left({\partial\over\partial y}\right)}^2\right)
$$
acting over the Hilbert space composed of all
$\Gamma$-automorphic functions which are square integrable with
respect to the hyperbolic measure ($\G = PSL(2,\,\ZZ))$.
Let $\{\psi_j\}_{j=1}^\infty$ be a maximal orthonormal system such that
$\Delta\psi_j=\lambda_j\psi_j$ for each $j\ge1$ and
$T(n)\psi_j=t_j(n)\psi_j$  for each integer $n\in\NN$, where
$$
\bigl(T(n)f\bigr)(z)\;=\;{1\over\sqrt{n}}\sum_{ad=n}
\,\sum_{b=1}^df\left({az+b\over d}\right)
$$
is the Hecke operator. We shall further assume that
$\psi_j(-\bar{z})=\e_j\psi_j(z)$ with $\e_j=\pm1$. We
then define ($s = \s + it$ will denote a complex variable)
$$
H_j(s)\;=\;\sum_{n=1}^\infty t_j(n)n^{-s}\qquad(\s > 1),
$$
which we call the Hecke series associated with the Maass wave form
$\psi_j(z)$, and which can be continued analytically
to an entire function over $\CC$.
It is known that $H_j(\hf) \ge 0$ (see Katok--Sarnak [15]), and that
$$
\sum_{\k_j\le K} \a_jH_j^3(\hf) = K^2\sum_{j=0}^3d_j\log^j K + O(K^{5/4}\log^{37/4}K)
\leqno(1.3)
$$
with suitable constants $d_j$, proved by the author in [9].
Here as usual  we insert in the sum over $\k_j$ the normalizing factor
$$
\a_j = |\rho_j(1)|^2(\cosh\pi\kappa_j)^{-1},
$$
where $\rho_j(1)$ is the first Fourier coefficient of  $\psi_j(z)$.
We also  have (see the author's paper [7])
$$
\sum_{K-G\le\k_j\le K+G} \a_j\H \;\ll_\e\; GK^{1+\e}\leqno(1.4)
$$
for
$$
K^{\e}  \;\le\; G \;\le \; K.\leqno(1.5)
$$
In view of $H_j(\hf) \ge 0$ we obtain from (1.4) the convexity-breaking bound
$H_j(\hf) \ll_\e \k_j^{1/3+\e}$, which is hitherto the sharpest one.

Note that by (1.3) and trivial estimation we obtain
$$
S(K) \;\ll\; K^2\log^3K,\leqno(1.6)
$$
and our wish is to try to  decrease the exponent of $K$ in (1.6).
It was conjectured in [8] that
$$
\sum_{K-1\le\k_j\le K+1} \a_j\H
\exp\left(i\k_j\log\left({\tau\over\k_j}\right)\right)
\;\ll_\e\; K^{1/2+\e}\leqno(1.7)
$$
holds for
$$
\tau^\delta \ll K \ll \tau^{1+\delta}\qquad(0 < \delta < 1).\leqno(1.8)
$$
This gives
$$
S(K) \;\ll_\e\; K^{3/2+\e},\leqno(1.9)
$$
thereby improving (1.6) by essentially a factor of $\sqrt{K}$. The conjecture
(1.7)--(1.8) is deep, and is certainly out of reach at present.
Heuristic reasons that it is best possible are given in [8].
It was also shown there
that its truth would imply essentially the best possible bounds for
the eighth moment of $|\zt|$, and for the error term (see (2.2)) in the
fourth moment formula for $|\zt|$.

\bigskip
\head
 2. Statement of results
\endhead
\bigskip
If $d(k)$ is the number of divisors of $k$, then we have

\medskip
THEOREM 1. {\it If $S(K)$ is defined by} (1.1) {\it and } (1.2)
{\it holds, then for some constants $0 < C_1 < C_2, c_\ell$
and $L\in\NN$, all of which
may be effectively evaluated, we have}
$$\eqalign{
S(K) &= \R\Biggl[\sum_{f\le3K}f^{{1\over2}}\sum_{C_1TK^{-1}f\le m
\le C_2TK^{-1}f}
m^{-{3\over2}}d(m)d(m+f){\roman e}^{i{Tf\over m}}\Bigl\{c_0 +\cr& +
\sum_{l=1}^Lc_\ell \f_\ell(K,T;m,f)\Bigr\}\Biggr] +
O_\e(K^{{3\over2}+\e}).\cr}\leqno(2.1)
$$
{\it The functions $\f_\ell(K,T;m,f)$ may be also explicitly evaluated,
and they are all $o(1)$ as $K\to\infty$ and} (1.2) {\it holds.}
\bigskip
The explicit shape of the functions $\f_\ell(K,T;m,f)$
will transpire during the proof, and a discussion on their precise shape
is given at the end of Section 5. Essentially they are (positive or
negative) powers in each variable. Thus
they are non-oscillating and, as stated,
all $o(1)$ as $K\to\infty$ and (1.2) holds. The
important fact is that they do not affect
the oscillating factor ${\roman e}^{i{Tf\over m}}$ in (2.1), and in fact can
be removed conveniently by partial summation techniques. For these reasons
it seemed more expedient to formulate Theorem 1 in the form given by (2.1),
than to write down explicitly all the functions $\f_\ell(K,T;m,f)$. The
number $L$ is a (large) constant, arising in (3.5) (and later in a similar
context). It comes from cutting the tails of a suitable series in such a way
that the tails in question make a negligible contribution.
By ``negligible contribution" we shall mean, here and later, a contribution
which is $\ll K_0^{-A}$ (or $\ll T^{-A}$) for any fixed $A>0$.

To abbreviate notation, sometimes in the proof we shall write expressions
similar to (2.1) as $A \asymp B + O_\e(K^{{3\over2}+\e})$. Namely $A\asymp B$
will mean, here and later, that $A$ is a multiple of $B$, plus a
finite number of sums (terms), each of which gives a bound not larger than
the bound for $B$, with some non-oscillatiing
functions $\f_\ell(K,T;m,f)$, as in (2.1).

\medskip
The importance of the sum $S(K)$ comes  primarily from its connection with
the function $E_2(T)$, the error term in the
asymptotic formula for the fourth moment of $|\zt|$.
This formula is customarily written as
$$
\int_0^T|\zt|^4\d t \;=\; TP_4(\log T) \;+\;E_2(T),
\quad
P_4(x) \;=\; \sum_{j=0}^4\,a_jx^j.\leqno(2.2)
$$
It was proved by A.E. Ingham that $a_4 = 1/(2\pi^2)$
(see e.g., [2, Chapter 5]), and much later by
D.R. Heath-Brown [1] that
$$
a_3 \;=\; 2(4\gamma - 1 - \log(2\pi) - 12\z'(2)\pi^{-2})\pi^{-2},
$$
who also produced
more complicated expressions for $a_0, a_1$ and $a_2$ in (2.3)
($\gamma = 0.577\ldots$ is Euler's constant). For an explicit evaluation of
the $a_j$'s  the reader is referred to [4].

In recent years, due  to the application of powerful methods
of spectral theory, much advance has been made in connection with
$E_2(T)$. We refer the reader to the works [3], [5], [6], [11]--[13], [20] and
[21]--[24]. Thus
N.I. Zavorotnyi [24] proved that $E_2(T) = O_\e(T^{2/3+\e})$, and it is
known now that
$$
E_2(T) \;=\; O(T^{2/3}\log^{C_1}T),\quad E_2(T) \;=\; \Omega(T^{1/2}),
\leqno(2.3) $$
and
$$
\int_0^TE_2(t)\d t \;=\; O(T^{3/2}),\quad\int_0^TE_2^2(t)\d t
\;=\; O(T^2\log^{C_2}T),\leqno(2.4)
$$
with effective constants $C_1,\,C_2 > 0$ (the values $C_1 = 8, C_2 = 22$ are
worked out in [23]). The above results were proved by
Y. Motohashi and the author (see [3], [11], [12] and [21]). The omega-result
in (2.3) ($f = \Omega(g)$ means that $f = o(g)$ does not hold,
$f = \Omega_\pm(g)$ means that $\limsup f/g > 0$ and that $\liminf f/g < 0$)
was improved to $E_2(T) = \Omega_\pm(T^{1/2})$ by Y. Motohashi [22].
There is an obvious discrepancy between the $O$--result and $\Omega$--result
in (2.3).  It was already mentioned that the conjecture
 $E_2(T) = O_\e(T^{1/2+\e})$
holds if the conjecture (1.7)--(1.8) is true.
It would imply  (by (2.9)) the hitherto unproved bound
$\ze \ll_\e t^{1/8+\e}$.

\medskip
Y. Motohashi proved (see [3, Chapter 6] and [23])
$$\eqalign{&
{1\over{\sqrt{\pi}}G}\int_{-\infty} ^{\infty}|\z(\hf+iT+it)|^4
\exp(-(t/G)^2)\d t\cr&
= {\pi\over\sqrt{2T}}\sum_{j=1}^\infty \a_j\H \k_j^{-{1\over2}}
\sin\left(\k_j\log{\k_j\over 4{\roman e}T}\right)\exp\Bigl(-{\txt{1\over4}}
\bigl({G\k_j\over T}\bigr)^2\Bigr) + O(\log^{3D+9}T),\cr}\leqno(2.5)
$$
if $\;T^{1/2}\log^{-D}T \le G \le T/\log T$ for an arbitrary, fixed
constant $D>0$, and
$$\eqalign{&
{1\over{\sqrt{\pi}}G}\int_0^V\int_{-\infty} ^{\infty}|\z(\hf+iT+it)|^4
\exp(-(t/G)^2)\d t\d T\cr&
= VP_4(\log V) +
\pi\sqrt{\hf V}\sum_{j=1}^\infty \a_j\H \k_j^{-{3\over2}}
\cos\left(\k_j\log{\k_j\over 4{\roman e}V}\right)\exp\Bigl(-{\txt{1\over4}}
\bigl({G\k_j\over V}\bigr)^2\Bigr) \cr&
+ O(V^{1/2}\log^CV) + O(G\log^5V),\cr}\leqno(2.6)
$$
for $\;V^{1/2}\log^{-A}V \le G \le V\exp(-\sqrt{\log V})$, $C = C(A)
\;(>0)$ for any  arbitrary, fixed constant $A>0$, where  $P_4$ is
given by (2.2). Then we have, as proved in [3, Lemma 5.1],
$$\eqalign{
E_2(2T) - E_2(T) &\le S(2T+\D\log T,\D) - S(T-\D\log T,\D)\cr&
+ O(\D\log^5T) + O(T^{1/2}\log^CT)\cr}\leqno(2.7)
$$
with $T^{1/2} \le\D\le T^{1-\e}$ and
$$
S(T,\D) := \pi\sqrt{\hf T}\sum_{j=1}^\infty \a_j\H \k_j^{-{3\over2}}
\cos\left(\k_j\log{\k_j\over 4{\roman e}T}\right)\exp\Bigl(-{\txt{1\over4}}
\bigl({\D\k_j\over T}\bigr)^2\Bigr).\leqno(2.8)
$$
A lower bound analogous to (2.7) holds also for $E_2(2T) - E_2(T)$,
and the estimation of $\zt$ is derived from [3, Lemma 4.1], namely
$$
\z(\hf + iT) \ll \log^{5/4}T + \left(\log T\max_{t\in[T-1,T+1]}
|E_2(t)|\right)^{1/4}.\leqno(2.9)
$$
The upper bound in (2.3) follows from (2.7)--(2.8) and trivial estimation,
namely (1.6), since the innocuous factors $\k_j^{-{3\over2}}$
and $\exp\Bigl(-{\txt{1\over4}}\bigl({\D\k_j\over T}\bigr)^2\Bigr)$
can be conveniently removed by partial summation from (2.8). Thus the
problem of the estimation of $E_2(T)$ (and hence also $\zt$) is reduced
to the estimation of our sum $S(K)$. The Lindel\"of exponent $\mu(\hf)$
is therefore seen not to exceed one fourth of the exponent in the bound
for $E_2(T)$ where, as usual, we define the Lindel\"of exponent as
$$
\mu(\s) \;=\;\limsup_{t\to\infty}\,{\log|\z(\s+it)|\over\log t }
\qquad(\s \in \RR).
\leqno(2.10)
$$
The famous, yet  unsettled, Lindel\"of hypothesis is that $\mu(\hf) =
0$, or equivalently that $\mu(\s) = 0$ for $\s\ge \hf$.

\medskip
The prominent feature of (2.1) is that the right-hand side contains
no quantities from spectral theory, but only classical exponential sums
with the divisor function $d(n) = \sum_{\delta|n}1$. In fact, the sum
in question can be considered as an exponential sum attached to the
so-called binary additive divisor problem (the evaluation and estimation
of $\sum_{m\le x}d(m)d(m+f)$, where $f$ is not fixed).
 Averages for $E(x;f)$, the error term in the asymptotic formula for this
sum, have been obtained by Y. Motohashi and the author  [13]. The
techniques developed in this work could be applied here,
since the problem reduces to the
evaluation of the sum ($X\approx Y$ means that $C_1X \le  Y \le C_2X$
holds for some constants $0 < C_1 < C_2$)
$$
\sum_{F<f\le2F}\int_N^{2N}{\roman e}^{i{Tf\over x}}E(x;f)\d x\qquad\left(F \ll K, \,
N \approx {TF\over K}\right).\leqno(2.11)
$$
Also the sum in (2.11) could be, at least in principle, evaluated
by Motohashi's formula [21] for the sum $\sum_{n=1}^\infty d(n)d(n+f)W(n/f)$,
where $W$ is a suitable smooth function.
Unfortunately, it appears that after the application of these procedures one
will eventually wind up with a sum of the same type as $S(K)$ in (1.1), plus some
manageable error terms.  The mechanism  is technically quite involved, and for
this reason it will not be discussed here in detail. However, it can be
seen heuristically from (4.4)--(4.7) of [13]. Namely the major contribution
to $E(x;f)$ comes from
$$
\R\left\{\hf x^{1/2}\sum_{\k_j\le Q}\a_j t_j(f)H_j^2(\hf)(f/x)^{i\k_j}
v(\k_j)\right\},\leqno(2.12)
$$
where $v(x) \ll x^{-3/2}$ and $Q$ is a  parameter satisfying certain
conditions. Inserting (2.12) expression in (2.11) we obtain exponential integrals
with the saddle point at $x_0 \approx TF/\k_j$, hence $\k_j \approx K$ is the relevant
range for $\k_j$. After the evaluation of the integral by the saddle point method
(see e.g., [2, Chapter 2]) we replace sums of $t_j(f)f^{-1/2}$ with $H_j(\hf)$
plus (small) error, to arrive at sums of the type $S(K)$ in (1.1), i.e., our
original sum.

\medskip
This type of impasse is well known from the estimation
of classical exponential sums (of the van der Corput type), where the
so-called $B$-process (essentially Poisson summation), when applied twice,
leads to the original exponential sum plus some (usually manageable)
error terms. It vitiates our efforts to attain a satisfactory estimate via
the application of binary additive problem techniques.
Naturally, one may  try other methods
to obtain from (2.1) a non-trivial bound, even if conditional estimates
such as the Lindel\"of hypothesis are assumed. However, at present
this seems difficult. One can separate the variables in (2.1) by setting
$n = m+f$ and letting $f$ lie in intervals of the form $[F,\,2F]$
with $F \ll K$. Then
the sum is majorized by $O(\log T)$ subsums of the form
$$
\left|\sum_{C_1TK^{-1}F\le m \le C_2TK^{-1}F}d(m)m^{-3/2}
\sum_{m+F<n\le m+2F}(n-m)^{1/2}d(n){\roman e}^{iTn/m}\right|.
$$
The factor $(n-m)^{1/2} $ can be conveniently removed by partial
summation. After that, one can apply the Vorono\"{\i} summation formula
(see e.g., [2, Chapter 3]) to the sum over $n$. The main difficulty
is that the sum over $n$ is ``short", in the sense that $F$ is much smaller
than $m$, and even after the application of the Vorono\"{\i} summation formula
to both sums, nothing better than the final trival estimate
$\ll_\e T^{1/2+\e}K^{3/2}$ seems to come out. This is no surprise, since
even the trivial bound
$$
\sum_{x<n\le x+h}d(n) \ll_\e hx^\e \;\qquad(1\ll h \ll x)
$$
cannot be obtained yet by the Vorono\"{\i} summation formula. Other methods,
such as the use of J.R. Wilton's approximate
functional equation and related
transformations (see M. Jutila [14]) can be also applied to the sum
over $n$, but the problem remains a very difficult one.

\bigskip
Instead of the sum $S(K)$ in (1.1) we may consider the analogous sums
when $\H$ is replaced by $H_j(\hf)$ or $H_j^2(\hf)$. The problem becomes
then considerably less difficult. On the other hand the exponential sums in
question do not seem to have immediate applications such as $S(K)$ does.
As we saw, $S(K)$ is crucial in the estimation of $E_2(T)$ and $\zt$,
which is our primary motivation. We shall prove

\bigskip
THEOREM 2. {\it If} (1.2) {\it holds, then}
$$
\eqalign{
\sum_{K<\k_j\le K'<2K} \a_j H_j^2(\hf)\cos\left(\k_j\log\left(
{4{\roman e}T\over \k_j}\right)\right) \;&\ll_\e \;T^{1/2+\e}K^{1/2},\cr
\sum_{K<\k_j\le K'<2K} \a_j H_j(\hf)\cos\left(\k_j\log\left(
{4{\roman e}T\over \k_j}\right)\right) \;&\ll_\e\; T^{1/2+\e}K^{1/4}.\cr}
\leqno(2.13)
$$

\bigskip
Therefore we see that the first bound improves the trivial bound
(see Y. Motohashi [23]) $O(K^2\log K)$ in the range
$T^{1/3+\e} \le K \le T^{1/2-\e}$. The  trivial bound for the
second sum in (2.13) is $O(K^2)$ (see Ivi\'c--Jutila [10]), and it
is improved for  $K$ satisfying $T^{2/7+\e} \le K \le T^{1/2-\e}$.
Clearly the method of proof of Theorem 1 and Theorem 2 can be used to
estimate certain other exponential sums of a similar nature.

\medskip
Similarly to the conjecture (1.9), one may conjecture that the sums
on the right-hand side of (2.13) are both $\ll_\e K^{3/2+\e}$. This conjecture,
like (1.9), is supported by the following mean square result. This is

\bigskip
THEOREM 3. {\it Let, for $m\in\NN$ and
$1 \ll K < K'\le 2K\ll T, T \le t \le 2T$,
$$
S_m(K;K',t) \;:=\;
\sum_{K<\k_j\le K'<2K} \a_j H_j^m(\hf)\cos\left(\k_j\log\left(
{4{\roman e}t\over \k_j}\right)\right).\leqno(2.14)
$$
Then, for $m = 1,2,3$,}
$$
\int_T^{2T}{\bigl(S_m(K;K',t)\bigr)}^2\d t \,\ll_\e\, T^{1+\e}K^3.\leqno(2.15)
$$

\bigskip
{\bf Corollary.} {\it We have}
$$
\int_0^T E_2^2(t)\d t \,\ll_\e\, T^{2+\e}.\leqno(2.16)
$$

\medskip
Note that (2.16) is a slightly weakened form of the second bound in (2.4),
obtained by Ivi\'c--Motohashi [11], and it is essentially
best possible (see the author's paper [6]).
The proof in [11] was based on a large
values estimate for $E_2(T)$, whose derivation
employed the spectral large sieve inequality. The new proof of (2.16) is simpler,
 being a direct consequence of (2.15) with $m=3$.

\medskip
The plan of the paper is as follows.
In Section 3 we make the technical preparation for the proof. Instead
of the ``long" sum (1.1), we shall use the transformation formulas involving
$H_j(\hf)$ for suitable (smooth) ``short" sums. Then we integrate over the
parameter to recover eventually the desired ``long" sum. The necessary
tool, which transforms our problem into a problem of the estimation
of the double exponential sum (cf. (2.1)) with two divisor functions,
is Motohashi's formula. It it presented in Section 4.
The principal part of the proof
of Theorem 1 is contained in Section 5, and the remaining part
will be given in Section 6.
Finally Theorem 2 is proved in Section 7, while Theorem 3 is proved in Section 8.
\head
 3. Technical preparation for the proof
\endhead
\bigskip
The basic idea of the proof of Theorem 1 is, as
with the proof of (1.4)-(1.5) in [7],
to use the transformation formula of Y. Motohashi (see [19] and
[23, Chapter Chapter 3]) for bilinear
forms of Hecke $L$-functions. Unfortunately, the shape (1.1)
of the fundamental sum $S(K)$ is not suited for the direct application
of the transformation formula. Before we can apply it, we have to
transform $S(K)$ into a suitable form. Although this is primarily a
technical problem, it is not obvious how one should tackle it,
and therefore the details will be given in this section.

\bigskip
We begin by considering, under the condition (1.2), the expression
$$\eqalign{
{1\over\sqrt{\pi}G}&\int_{K_0}^{K_0'}
\sum_{j=1}^\infty \a_j\H \exp\left(i\k_j\log{4{\roman e}T\over\k_j}
-(\k_j-K)^2G^{-2}\right)\d K\cr&
= {1\over\sqrt{\pi}G}\int_{K_0}^{K_0'}\sum(K;G)\d K,\cr}\leqno(3.1)
$$
say, where $G = G(K_0)$ is a parameter satisfying
$$
K_0^\e \le G \le K_0^{1/2-\e},\quad K_0 \le K \le K_0' \le 2K_0.\leqno(3.2)
$$
Exchanging the order of summation and integration in (3.1) we have, in view of (1.1),
that
$$\eqalign{&
{1\over\sqrt{\pi}G}\R\left\{\int_{K_0}^{K_0'}\sum(K;G)\d K\right\}\cr&
= {1\over\sqrt{\pi}G}\R\Biggl\{
\sum_{K_0-G\log K_0\le \k_j\le K_0'+G\log K_0}
\a_j\H\exp\left(i\k_j\log{4{\roman e}T\over \k_j}\right)\times\cr&\times
\int_{\k_j-G\log K_0}^{\k_j+G\log K_0}
{\roman e}^{-(\k_j-K)^2G^{-2}} \d K\Biggr\} + O_\e(K_0^\e)\cr&
= \R\Biggl\{\sum_{K_0-G\log K_0\le \k_j\le K_0'+G\log K_0}
\a_j\H\exp\left(i\k_j\log{4{\roman e}T\over \k_j}\right)\times\cr&\times
{1\over\sqrt{\pi}}\int_{-\log K_0}^{\log K_0}{\roman e}^{-u^2}\d u\Biggr\}
+ O_\e(K_0^\e)\cr&
= S(K_0;K_0',T) + O_\e(K_0^{1+\e}G),\cr}\leqno(3.3)
$$
where we used (1.4) to estimate the contribution from $\k_j$ lying
in the intervals $[K_0-G\log K_0,\, K_0]$ and $[K_0',\, K_0'+G\log K_0]$.
On the other hand we have
$$\eqalign{&
{1\over\sqrt{\pi}G}\int\limits_{K_0}^{K_0'}\sum(K;G)\d K = O_\e(K_0^\e)
\,+\cr&
+  {1\over\sqrt{\pi}G}\int\limits_{K_0}^{K_0'}
\sum_{|\k_j-K|\le G\log K_0}
\a_j\H\exp\left(i\k_j\log{4{\roman e}T\over \k_j}
-(\k_j-K)^2G^{-2}\right) \d K.\cr}
\leqno(3.4)
$$
Note that, for $|\k_j-K|\le G\log K_0$ and $K_0 \le K \le 2K_0$, we have
$$\eqalign{&
\k_j\log{4{\roman e}T\over \k_j} - K - \k_j\log{4T\over K}
= \k_j - K + \k_j\log {K\over\k_j}\cr&
= \k_j - K + \k_j\log\left(1 +  {K-\k_j\over\k_j}\right)\cr&
= \sum_{\ell=2}^L\,{(-1)^{\ell-1}\over \ell}\k_j
\left({K-\k_j\over\k_j}\right)^{\ell}
+ O\left({G^{L+1}\log^{L+1}K_0\over K_0^L}\right)\cr}\leqno(3.5)
$$
for any fixed integer $L\ge 2$. But as, for $\ell\in\NN$,
$|\k_j-K|\le G\log K_0$,
$$
\k_j^{-\ell} = K^{-\ell}\left(1 + {\k_j-K\over K}\right)^{-\ell}
= K^{-\ell}\left\{1 + \sum_{j=1}^\infty {-\ell\choose j}
\left({\k_j-K\over K}\right)^\ell\right\},
$$
we obtain
$$
\eqalign{&
\exp\Bigl(i\k_j\log{4{\roman e}T\over \k_j} \Bigr)
=\exp\Biggl\{iK + i\k_j\log{4T\over K}\cr&
+  i\sum_{\ell=2}^L\,{(-1)^{\ell-1}\over \ell}\k_j\Bigl({K-\k_j
\over\k_j}\Bigr)^{\ell}
+ O\left({G^{L+1}\log^{L+1}K_0\over K_0^L}\right)\Biggr\}\cr&
= {\roman e}^{iK}\exp\left(i\k_j\log{4T\over K}\right)\cdot
\Biggl\{1 + \sum_{\ell=2}^L a_\ell K^{1-\ell}(K-\k_j)^{\ell}+\cr&
+ O\left({G^{L+1}\log^{L+1}K_0\over K_0^L}\right)\Biggr\}
\cr}\leqno(3.6)
$$
with suitable constants $a_\ell$. In view of (3.2) we can choose
$L\,(\ge 2)$ so large that the error term in (3.6), when inserted in
(3.4), will make a contribution which is negligible (i.e., $\ll K_0^{-A}$
for any given $A>0$).

\medskip
The remaining terms in (3.6) have the property that the summands in
the sum over $\ell$ are of decreasing order of magnitude, since
for $|K  -\k_j| \le G\log K_0$ and  $K_0 \le K\le2K_0$, we have
$$
|K  -\k_j|K^{-1} \ll GK_0^{-1}\log K_0 \ll_\e K_0^{-\e-1/2}\log K_0.
$$
Therefore we can write
$$
\eqalign{&
{1\over\sqrt{\pi}G}\int\limits_{K_0}^{K_0'}
\sum_{|\k_j-K|\le G\log K_0}
\a_j\H{\roman e}^{iK}\exp\left(i\k_j\log{4T\over K}\right)
{\roman e}^{-(\k_j-K)^2G^{-2}} \d K\cr&
= {1\over\sqrt{\pi}G}\int\limits_{K_0}^{K_0'}
 R_0(K;T,G) {\roman e}^{iK}\cdot\d K \cr&
 + \sum_{\ell=2}^La_\ell{1\over\sqrt{\pi}G}\int\limits_{K_0}^{K_0'}
 R_\ell(K;T,G) {\roman e}^{iK}\cdot\d K+ O_\e(K_0^\e),\cr}\leqno(3.7)
$$
say, where for $\ell = 0,1,2,\ldots$ we have set
$$
R_\ell(K;T,G) := \sum_{j=1}^\infty\a_j\H h_\ell(\k_j;T,K,G),\leqno(3.8)
$$
and the function $h_\ell$ is defined as follows. For a fixed $N\in\NN$ we set
$$
q_N(r) := {\left(r^2+{1\over4}\right)\left(r^2+{9\over4}\right)\cdots
\left(r^{2}+ {(2N-1)^2\over4}\right)\over (r^2+ 100N^2)^N},
\leqno(3.9)
$$
and then define
$$\eqalign{
h_\ell(r;T,K,G) &:= q_N(r)\left(L_\ell(r;T,K,G) + L_\ell(-r;T,K,G)\right),\cr
L_\ell(r;T,K,G) &:= K^{1-\ell}(K-r)^\ell\left({4T\over K}\right)^{ir}
{\roman e}^{-(r-K)^2G^{-2}},\cr}\leqno(3.10)
$$
so that $h_\ell$ is an even function of $r$.
From (3.3) and (3.5)--(3.10) it follows that
$$\eqalign{&
S(K_0;K_0',T) =
\R\left\{{1\over\sqrt{\pi}G}\int\limits_{K_0}^{K_0'}
 R_0(K;T,G) {\roman e}^{iK}\cdot\d K \right\}\cr&
 + \R\left\{\sum_{\ell=2}^La_\ell{1\over\sqrt{\pi}G}\int\limits_{K_0}^{K_0'}
 R_\ell(K;T,G) {\roman e}^{iK}\cdot\d K\right\}+ O_\e(K_0^{1+\e}G),\cr}\leqno(3.11)
$$
and clearly the main contribution to our sum (cf. (1.1))
 $S(K_0;K_0',T)$ comes from the integral with $R_0$.

\medskip
The function $h_\ell(r;T,K,G)$, defined by (3.10),
is a modified Gaussian weight function in $r$, which is regular in the
horizontal strip $|\I r| \le N+1$. Moreover it is even,
satisfies $h_\ell(\pm\hf ij;T,K,G) = 0$
for $j = 1,3,\ldots, \hf(N-1)$ and every $\ell$ and the decay condition
$$
h_\ell(r;T,K,G) \;\ll_{\ell,T,K,G}\; \exp(-c|r|^2)\qquad(c>0)\leqno(3.12)
$$
in the above strip.  Thus it satisfies all the conditions necessary
for the application of Motohashi's transformation formula, which will
be discussed in the next section. This ends the technical preparation
for the proofs.

\bigskip
\head
 4. Motohashi's transformation formula
\endhead
\bigskip
The basic idea of the transformation formula is to transform
the expression, for a suitable weight function $h_0(r)$,
$$
{\Cal C}(K,G) := \sum_{j=1}^\infty \a_j\H h_0(\k_j)        \leqno(4.1)
$$
into a sum of terms which do not contain quantities from the
spectral theory of the non-Euclidean Laplacian. In this way the problem
of the evaluation or estimation of ${\Cal C}(K,G)$ is transformed into
a problem of classical Analytic Number Theory. The
function ${\Cal C}(K,G)$
will be actually $R_\ell(K;T,G)$ from (3.8). For the function $h_0(r)$,
which is regular in a (large) fixed horizontal strip, it is sufficient
to assume that it is even and decays in the strip like
$$
h_0(r) \;\ll\; {\roman e}^{-c|r|^2}\qquad(c>0).\leqno(4.2)
$$
We set  $\lambda = C\log K\;(C > 0)$ and note that one has
(this is Y. Motohashi [23, eq. (3.4.18)], with the
extraneous factor $(1 - (\k_j/K)^2)^\nu$ omitted)
$$\eqalign{
{\Cal C}(K,G) &= \sum_{f\le3K}f^{-{1\over2}}\exp\Bigl(-{\bigl({f\over K}
\bigr)}^\lambda\Bigr){\Cal H}(f;h_0)\cr&
- \sum_{\nu=0}^{N_1}\,\sum_{f\le3K}f^{-{1\over2}}U_\nu(fK)
{\Cal H}(f;h_\nu) + O(1),\cr}                    \leqno(4.3)
$$
with
$$\eqalign{
h_\nu(r) &\= h_0(r)\left(1 - \left({r\over K}\right)^2\right)^\nu,\cr
{\Cal H}(f;h) &\= \sum_{\nu=1}^7{\Cal H}_\nu(f;h),\cr}\leqno(4.4)
$$
$$
{\Cal H}_1(f;h) \= -2\pi^{-3}i\left\{(\gamma - \log(2\pi{\sqrt f}))
({\hat h})'(\hf) + {\txt{1\over4}}({\hat h})''(\hf)\right\}d(f)
f^{-{1\over2}},
$$
$$
{\Cal H}_2(f;h) \= \pi^{-3}\sum_{m=1}^\infty m^{-{1\over2}}d(m)d(m+f)
\Psi^+\Bigl({m\over f};h\Bigr),
$$
$$
{\Cal H}_3(f;h) \= \pi^{-3}\sum_{m=1}^\infty (m+f)^{-{1\over2}}
d(m)d(m+f)\Psi^-\Bigl(1+ {m\over f};h\Bigr),        \leqno(4.5)
$$
$$
{\Cal H}_4(f;h) \= \pi^{-3}\sum_{m=1}^{f-1}m^{-{1\over2}}d(m)d(f-m)
\Psi^-\Bigl({m\over f};h\Bigr),
$$
$$
{\Cal H}_5(f;h) \= -(2\pi^3)^{-1}f^{-{1\over2}}d(f)\Psi^-(1;h),
$$
$$
{\Cal H}_6(f;h) \= -12\pi^{-2}i\s_{-1}(f)f^{1\over2}h'(-\hf i),
$$
$$
{\Cal H}_7(f;h) \= -\pi^{-1}\int_{-\infty}^\infty{|\z(\hf+ir)|^4\over
|\z(1+2ir)|^2}\s_{2ir}(f)f^{-ir}h(r)\d r\quad (\s_a(f) = \sum_{d|f}d^a),
$$
where
$$
{\hat h}(s) \= \int_{-\infty}^\infty rh(r){\G(s+ir)\over\G(1-s+ir)}\d r,
\leqno(4.6)
$$
$$
\Psi^+(x;h) \= \int_{(\b)}\G^2(\hf-s)\tan(\pi s){\hat h}(s)x^s\d s,
\leqno(4.7)
$$
$\int_{(\b)}$ denotes integration over the line $\R s = \b$,
$$
\Psi^-(x;h) \= \int_{(\b)}\G^2(\hf-s){{\hat h}(s)\over
\cos(\pi s)}x^s\d s\leqno(4.8)
$$
with $-{3\over2} < \b < {\hf}$, $N_1$ is a sufficiently large integer,
$$
U_\nu(x) = {1\over2\pi i\lambda}\int_{(-\lambda^{-1})}
(4\pi^2K^{-2}x)^wu_\nu(w)
\G({w\over\lambda})\d w \ll \left({x\over K^2}\right)^{-{C\over\log K}}
\log^2K,
$$
where $u_\nu(w)$ is a polynomial in $w$ of degree $\le 2N_1$, whose
coefficients are bounded. As already mentioned,
the  prominent feature of Motohashi's explicit
expression for ${\Cal C}(K,G)$ is that it contains series and integrals
with the classical divisor function $d(n)$ only, with no quantities
from spectral theory.

\bigskip
\head
 5. Proof of Theorem 1
\endhead
\bigskip
We need, in view of (3.11), to transform and estimate
the functions $R_\ell(K;T,G)$ in (3.8). To this end we shall
employ (4.1), where $h_0(r)$ equals
$$
h_\ell(r;T,K,G){\left(1 - \left({r\over K}
\right)^2\right)}^\nu\quad(\nu = 0,1,2,\ldots\,;
\ell = 0,1,2,\ldots\,).\leqno(5.1)
$$
All the functions of the form (5.1) are treated analogously.
Therefore it is sufficient to consider in detail only
the case $\nu = \ell = 0$, when  for simplicity the function in
(5.1) will be again denoted by $h(r)$. It is clearly this case
which will furnish eventually the largest contribution to (2.1).

\smallskip
In the sequel we shall repeatedly use the classical formula
$$
\int_{-\infty}^\infty {\roman e}^{Au-Bu^2}\d u
\= \sqrt{\pi\over B}\exp\left({A^2\over4B}\right)
\qquad(\R B > 0).\leqno(5.2)
$$
By taking $B=1$ and then differentiating (5.2) as the function of $A$,
we also obtain
$$
\int_{-\infty}^\infty u^j{\roman e}^{Au-u^2}\d u
\= P_j(A){\roman e}^{{1\over4}A^2}
\qquad(j = 0,1,2,\ldots,\;P_0(A) = \sqrt{\pi}\,),\leqno(5.3)
$$
where $P_j(z)$ is a polynomial in $z$ of degree $j$, which may be
explicitly evaluated. The basic idea is that
the factor $(4T/K)^{\pm ir}$ (cf. (3.10))
is the dominating oscillating factor which in most cases,
after the use of (5.2)
or (5.3), will produce exponential functions of fast decay which
will make a negligible contribution. We recall that a
``negligible contribution" is one which is $\ll K_0^{-A}$ (or
$\ll T^{-A}$) for any fixed $A>0$.

\medskip
This is precisely what happens with  the contribution of ${\Cal H}_1(f;h)$,
which we shall first show to be negligible. Namely
from (4.6) we find that
$$
(\hat h)'(\hf) = 2\int_{-\infty}^\infty rh(r){\G'\over\G}(\hf+ir)\d r.
\leqno(5.4)
$$
But (see e.g., [18])
$$
{\G'\over\G}(s) = \log s - {1\over2s} + O\left({1\over|s|}\right),
\leqno(5.5)
$$
where the $O$-term admits an asymptotic expansion. The non-negligible
contribution in (5.4) is for the range $|r\pm K|\le G\log K$. We
make the change of variable $r \pm K = Gu$ and use Taylor's formula
to simplify the integrand. After this we may use (5.2) and (5.3), which
will produce exponential factors of the form $\exp(-{1\over4}G^2
(\log{4T\over K})^2)$, which will make a negligible
contribution. The $O$-term in (5.5), by trivial estimation, will
make a total contribution of $O_\e(K^{3/2+\e})$. The
contribution of $({\hat h})''(\hf)$ is estimated analogously, and we see
that the total contribution of ${\Cal H}_1(f;h)$ is $O_\e(K_0^{3/2+\e})$.

Next we note that
$$
{\Cal H}_6(f;h) \= -12\pi^{-2}i\s_{-1}(f)f^{1\over2}h'(-\hf i)
\ll \s_{-1}(f)f^{1\over2}\exp(-\hf K^{2}G^{-2}),
$$
hence summation over $f$ in (4.3) yields a contribution which is
negligible.

The total contribution of
$$
{\Cal H}_5(f;h) \= -(2\pi^3)^{-1}f^{-{1\over2}}d(f)\Psi^-(1;h)\leqno(5.6)
$$
is also negligible. This follows from [23, eq. (3.3.44)], in view of the
presence of $\sinh\pi r/(\cosh\pi r)^2$, which decays like $\exp(-\pi |r|)$.

The total contribution of
$$
{\Cal H}_3(f;h) \= \pi^{-3}\sum_{m=1}^\infty (m+f)^{-{1\over2}}
d(m)d(m+f)\Psi^-\left(1+ {m\over f};h\right)  \leqno(5.7)
$$
is also negligible, but this is somewhat more involved than the
contribution of ${\Cal H}_5(f;h)$.  We need the representation (this
is [23, eq. (3.3.43)])
$$
\Psi^-(x;h) = 2\pi i\int_0^1(y(1-y)(1-y/x))^{-1/2}\int_{-\infty}^\infty
{rh(r)\over\cosh(\pi r)}\left\{{y(1-y)\over x-y}\right\}^{ir}
\d r\d y,\leqno(5.8)
$$
which is valid for $x>1$. Motohashi derived (5.8) for a somewhat different
weight function $h(r)$, essentially without the factor
$(4T/K)^{\pm ir}$, but it is clear
by following his proof that (5.8) will hold for the present function $h(r)$
as well. The same remark holds for other forms of the functions
$\Psi^\pm(x;h)$ which will be needed in the sequel.
To deal with the series over $m$ in (5.7) we need to have a good
bound in $m$. This is achieved, as in [23], by shifting the line of
integration (in the integral over $r$) in (5.8) to $\I r = -1$. In this
process use is made of the fact that $h(-\hf i) = 0$, since this zero
at $-\hf i$ cancels with the zero of $\cosh\pi r$.
We then note that, in the relevant range for $r$,
$1/\cosh(\pi r) \ll \exp(-\hf\pi K)$. Thus, for $x = 1+ m/f \ge 3$, we
obtain by trivial estimation
$$
\Psi^-\left(1 + {m\over f};h\right) \ll fm^{-1}TG\exp(-\hf\pi K)
\qquad(m\ge 2f).
$$
This is more than sufficient to render the total contribution of $m \ge 2f$
negligible, and the same follows for the contribution of the remaining
$m$'s if we use the trivial estimate (coming directly from (5.8))
$$
\Psi^-\left(1 + {m\over f};h\right) \ll KG\exp(-\hf\pi K)\qquad(m\le 2f).
$$

\medskip
To deal with
$$
{\Cal H}_7(f;h) \= -\pi^{-1}\int_{-\infty}^\infty{|\z(\hf+ir)|^4\over
|\z(1+2ir)|^2}\s_{2ir}(f)f^{-ir}h(r)\d r,
$$
note that we have
 $1/\z(1+ir) \ll \log (|r| + 1)$ and
$$
\sum_{n=1}^\infty \s_{2ir}(n)n^{-ir-s} \= \z(s-ir)\z(s+ir) \qquad
(r \in \RR,\, \R s > 1).
$$
Consequently by the Perron inversion formula (see e.g., [2, eq. (A.10)])
$$
\sum_{f\le3K}\,\s_{2ir}(f)f^{-{1\over2}-ir} \;\ll_\e\;
K^{2\mu({1\over2})+\e} \;\ll_\e\; K^{{1\over3}+\e}
\qquad (K \ll |r| \ll K),\leqno(5.9)
$$
where $\mu(\s)$  is given by (2.10),
and we used the classical bound $\mu(\hf) \le 1/6$.
Since the relevant range of $r$ in ${\Cal H}_7(f;h)$ is
$|r \pm K| \le G\log K_0$, it follows by using (5.9) that
$$\eqalign{&
G^{-1}\int_{K_0}^{K'_0}\sum_{f\le3K}f^{-1/2}{\Cal H}_7(f;h)\d K
\cr&\ll_\e 1+
K_0^{1/3+\e}G^{-1}\int_{K_0}^{K'_0}\int_{K-G\log K_0}^{K+G\log K_0}|\z(\hf+ir)|^4
\d r\,\d K \cr&
\ll_\e K_0^{1/3+\e}G^{-1}\int_{K_0-G\log K_0}^{K'_0+G\log K_0}|\z(\hf+ir)|^4
\int_{r-G\log K_0}^{r+G\log K_0}\d K\cdot\d r\cr&
\ll_\e K_0^{4/3+\e},\cr}
$$
hence this is the total contribution of ${\Cal H}_7(f;h)$
to the right-hand side of (3.7).

\medskip
It remains yet to deal with
$$
{\Cal H}_2(f;h) \= \pi^{-3}\sum_{m=1}^\infty m^{-{1\over2}}d(m)d(m+f)
\Psi^+\left({m\over f};h\right)\leqno(5.10)
$$
and ${\Cal H}_4(f;h)$, which will be done in Section 6.
The contribution of ${\Cal H}_2(f;h)$ is
the principal one. It is estimated
according to the range of $m$ in (5.10).

\medskip
We shall show first that the contribution of $m \ge fTK^{\e-1}$ in
(5.10) is negligible. We use the representation (this is [23, eq. (3.3.41)])
$$
\Psi^+(x;h) = 2\pi\int_0^1\{y(1-y)(1+y/x)\}^{-1/2}
\int_{-\infty}^\infty rh(r)\tanh(\pi r)\left\{{y(1-y)\over x+y}
\right\}^{ir}\d  r\d y\leqno(5.11)
$$
with $x = m/f \ge K^\e$, and shift the line of integration in the
inner integral to $\I r = -N$. This is permissible, since by (3.9)
and (3.10) the function $h(r)$ is regular for $|\I r| \le N+1$. Then
the inner  integral in (5.11) becomes
$$\eqalign{&
\int_{-\infty}^\infty (r-Ni)h(r-Ni)\tanh(\pi r)\left\{{y(1-y)\over x+y}
\right\}^{ir}\left\{{y(1-y)\over x+y}\right\}^N\d r\cr&
\ll KG(y(1-y))^N\left({Tf\over mK}\right)^N.\cr}
$$
Since $N \;(= N(\e))$ can be taken arbitrarily large, it follows that
the total contribution of $m/f \ge TK^{\e-1}$ in (5.10) is negligible.

We shall show that the contribution of $m/f \le TK^{-\e-1}$
is also negligible. We make the change of variable $r = \pm K + Gu$ in the
$r$-integral in (5.11), and note that
$$
\tanh(\pi r) \= \roman {sgn}(r) + O({\roman e}^{-2\pi|r|})
\qquad(r\in\RR).
\leqno(5.12)
$$
After the application of (5.2)
there will appear the exponential factors
$$
\exp\left(-{\txt{1\over4}}G^2\log^2\left({4T\over K}\cdot
{y(1-y)\over x + y}\right)\right)
$$
and
$$
\exp\left(-{\txt{1\over4}}G^2\log^2\left({4T\over K}\cdot
{ x + y\over y(1-y)}\right)\right).
$$
Since, in view of (1.2),
$$
{4T\over K}\cdot
{ x + y\over y(1-y)} \ge {4Tx\over K} = {4Tm\over fK} \ge {4T\over 3K^2}
\gg T^{\e/2},
$$
the contribution of the latter is negligible. The contribution of the
former is also negligible if
$$
{4T\over K}\cdot
{y(1-y)\over x+y} \le 1 - G^{-1}\log T\quad{\roman {or}}\quad
{4T\over K}\cdot
{  y(1-y)\over x+y} \ge 1 + G^{-1}\log T.
$$
If this condition is not satisfied, then
$$
y \in [y_1,\,y_2],\; y_1 \approx Kx/T \ll K^{-\e},\;
y_1 - y_2 \approx {Kx\log T\over TG}.
$$
In the $y$-integral in (5.11) over $[y_1,\,y_2]$
we integrate by parts the factor
$y^{ir-{1\over2}}$ a large number of times. Each time the exponent
of $y$ will increase by unity, while the order of the $r$-integral will
remain unchanged. Trivial estimation shows then that the contribution
of $m/f \le TK^{-\e-1}$ is indeed negligible.
\medskip
Thus the critical range in the estimation of ${\Cal H}_2(f;h) $ is
(since $K_0 \le K \le 2K_0$)
$$
fTK_0^{-1-\e} \;\le m \le \;fTK_0^{-1+\e}.\leqno(5.13)
$$
For the range (5.13) we shall use the representation which
follows from [23, eq. (3.3.39)] and the formula after it,
with $x = m/f \to\infty$ (as $K_0\to\infty$), namely
$$\eqalign{&
\Psi^+(x;h) =\cr& 2\pi\int_{-\infty}^\infty rh(r)\tanh(\pi r)
\R\left\{{\G^2(\hf + ir)\over\G(1+2ir)}
F\left(\hf+ir,\hf+ir;1+2ir;-{1\over x}\right)x^{-ir}\right\}\d r,\cr}
\leqno(5.14)
$$
where $F$ is the hypergeometric function. We could use the asymptotic
formula, valid for $y \ge y_0 > 1$ and $r\to\infty$,
$$\eqalign{&
F\left(\hf+ir,\hf+ir;1+2ir;-{1\over y^2}\right) = O(y^{-4}r^{-2})\cr&
+ (2y)^{2ir}(y + \sqrt{1+y^2})^{-2ir}\left({y^2\over1+y^2}\right)^{1/4}
\left(1 - {1\over8ir}\cdot{2y^2+1\over2y\sqrt{1+y^2}}\right),\cr}
\leqno(5.15)
$$
which yields directly the main term.
This formula is to be found in the work of N.I. Zavorotnyi [24]. A sketch
of proof is given by N.V. Kuznetsov [17], where the asymptotics are given
by means of a solution of a certain second-order differential equation
(see his work [16]).
One can avoid the use of (5.15) by appealing to the  classical
quadratic transformation formula (see [18, eq.
(9.6.12)]) for the hypergeometric
function, as was done by the author [7] in his work on sums of $\a_j\H$
in short intervals. This is
$$
F(\a,\b;2\b;z) \= \left({1+\sqrt{1-z}\over2}\,\right)^{-2\a}
F\left(\a,\a-\b+\hf;\b+\hf;\left({1-\sqrt{1-z}\over1+\sqrt{1-z}}\,\right)^2
\right),
$$
and then one can develop the resulting hypergeometric function into
a convergent power series, of which the main contribution will come
from the leading term, namely 1. The main term in (2.1) (the summand with
$c_0$) will be in both cases the same, of course, and the latter approach
yields the remaining summands with $\f_\ell$.

\medskip
In (5.14) the relevant ranges of integration are $[-K - G\log K_0,
-K + G\log K_0]$ and $[K - G\log K_0,
K + G\log K_0]$. We note that

and in the first range of integration we change $r$ to $-r$. Then we obtain
 that the critical expression in question is
$$\eqalign{&
{4\sqrt{\pi}\over G}\int_{K_0}^{K_0'}{\roman e}^{iK}
\sum_{f\le3K_0}f^{-1/2}\sum_{TK_0^{-1-\e}f\le m \le TK_0^{-1+\e}f}
m^{-1/2}d(m)d(m+f)\times\cr&
\int_{K-G\log K_0}^{K+G\log K_0}r{\left({4T\over K}\right)}^{ir}
{\roman e}^{-(r-K)^2G^{-2}}\times\cr&
\R\left\{{\G^2(\hf + ir)\over\G(1+2ir)}
F\left(\hf+ir,\hf+ir;1+2ir;-{1\over x}\right)
x^{-ir}\right\}\d r\d K.\cr}\leqno(5.16)
$$
To (5.16) we shall apply (5.15) with $y = \sqrt{x} = \sqrt{m/f}$,
under (5.13). The gamma-factors are simplified by Stirling's
formula, namely that for $t \ge t_0 >0$
$$
\G(s) = \sqrt{2\pi}\,t^{\s-{1\over2}}\exp\left(-\hf\pi t + it\log t -it
+ \hf{\pi i}(\s - \hf)\right)\cdot\left(1 + O_\s\left(
t^{-1}\right)\right),\leqno(5.17)
$$
with the understanding that the $O$--term in (5.17) admits an
asymptotic expansion in terms of negative powers of $t$. Hence using
the symbol $\asymp$ (defined after the formulation of Theorem 1)
the expression in (5.16) is ($x = m/f)$
$$\eqalign{&
\asymp {1\over G}\int_{K_0}^{K_0'}{\roman e}^{iK}\sum_f\sum_m\cdots
\int_{K-G\log K_0}^{K+G\log K_0}r{\left({4T\over K}\right)}^{ir}
{\roman e}^{-(r-K)^2G^{-2}}\times\cr&
\R\left\{r^{-1/2}{\roman e}^{-2ir\log 2}x^{-ir}2^{2ir}x^{ir}
(\sqrt{x} + \sqrt{1+x})^{-2ir}\d r\right\}\d K
\cr&
\asymp {1\over G}\int_{K_0}^{K_0'}K^{1/2}{\roman e}^{iK}
\sum_f\sum_m\cdots\times
\cr& \times\int_{K-G\log K_0}^{K+G\log K_0}{\left({4T\over K}\right)}^{ir}
\cos(2r\log(\sqrt{x} + \sqrt{1+x}))
{\roman e}^{-(r-K)^2G^{-2}}\d r\d K.
\cr}\leqno(5.18)
$$
The cosine is written as the sum of exponentials, after which the
change of variable $r = K + Gu$ is made in the $r$-integral.
The inner integral in (5.18) thus reduces to
$$
G\int_{-\log K_0}^{\log K_0}{\roman e}^{-u^2}
\exp\left\{(iK+iGu)\left(\log{4T\over K}\pm
\log(\sqrt{x} + \sqrt{1+x})^2\right)\right\}\,\d u,\leqno(5.19)
$$
after which we restore the integration to the whole real line,
making a negligible error. Then we apply (5.2), noting that the
integral with the $+$-sign makes a negligible contribution.
The integral with the  $-$-sign equals
$$
\sqrt{\pi}G\exp\left\{iK\log\left({4T\over K(\sqrt{x} + \sqrt{1+x})^2}
\right) - {1\over4}G^2\log^2\left({4T\over K
(\sqrt{x} + \sqrt{1+x})^2}\right)\right\}.
$$

It follows that (5.18) is
$$\eqalign{&
\asymp \sum_f\sum_m\cdots\int_{K_0}^{K_0'}K^{1/2}
\exp\left\{iK\log\Bigl({4{\roman e}T
\over K(\sqrt{x} + \sqrt{1+x})^2}\Bigr)\right\}
\cr&\times
\exp\left\{-{\txt{1\over4}}G^2
\log^2\Bigl({4T\over K(\sqrt{x} + \sqrt{1+x})^2}\Bigr)\right\}\d K.
\cr}\leqno(5.20)
$$
The last exponential factor yields that only the contribution of
$m/f \approx T/K_0$ makes a non-negligible contribution.
More precisely, we have
$$
{4T\over K(\sqrt{x} + \sqrt{1+x})^2}
= {T\over K\left(x + \sum\limits_{j=0}^\infty b_jx^{-j}\right)}
\qquad(x = m/f > 1)
$$
with suitable coefficients $b_j$. Therefore the second exponential
factor in (5.20) is negligibly small, unless
$$
K = {T\over x + \sum\limits_{j=0}^\infty b_jx^{-j}}\left(1 +
O\left({\log T\over G}\right)\right).\leqno(5.21)
$$
This means that the relevant interval of integration over $K$
in (5.20), for fixed $f$ and $m$, has length $\ll Tf\log T/(mG)$.

The integral in (5.20) is an exponential integral of the form
$$\eqalign{&
\int_{K_0}^{K_0'} g(K){\roman e}^{if(K)}\d K,\quad
f(K) := K\log\Bigl({4{\roman e}T
\over K(\sqrt{x} + \sqrt{1+x})^2}\Bigr).\cr&
g(K) := K^{1/2}\exp\left\{-{\txt{1\over4}}G^2
\log^2\Bigl({4T\over K(\sqrt{x} + \sqrt{1+x})^2}\Bigr)\right\}.
\cr}
$$
The saddle point $K_1$ (the root of $f'(K)=0)$ is given by
$$
K_1 = {4T\over (\sqrt{x} + \sqrt{1+x})^2}.\leqno(5.22)
$$
Since $f''(K) = -1/K$, it follows by the saddle
point method (see e.g., [2, Chapter 2])
 that (5.20) is ($0 < C_1 < C_2$ are suitable constants, $x = m/f$)
$$
\asymp T
\sum_{f\le3K_0}f^{-\hf}\sum_{{C_1Tf\over K_0}\le m \le {C_2Tf\over K_0}}
m^{-\hf}{d(m)d(m+f)\over(\sqrt{x} + \sqrt{1+x})^2}
\exp\left({4iT\over (\sqrt{x} + \sqrt{1+x})^2}
\right),
$$
plus an error term which is certainly $\ll_\e K_0^{3/2+\e}$. But since
$$
{4iT\over (\sqrt{x} + \sqrt{1+x})^2}
= {iT\over x}\left(1 + \sum_{j=1}^\infty c_j x^{-j}\right)\leqno(5.23)
$$
with suitable constants $c_j$ and $Tx^{-2} \ll K^2/T \ll T^{-\e}$ in
view of (1.2), it follows that (5.20) is
$$
\asymp T
\sum_{f\le3K_0}f^{1\over2}\sum_{{C_1Tf\over K_0}\le m \le {C_2Tf\over K_0}}
m^{-{3\over2}}d(m)d(m+f)\exp\left({iTf\over m}
\right) + O_\e(K_0^{3/2+\e}).\leqno(5.24)
$$
Therefore the proof of Theorem 1 will be complete after we show that the
contribution of ${\Cal H}_4(f;h) $ is negligible,
and choose $G = K_0^{1/2-\e}$. Note that trivial estimation gives
that the expression in (5.24) is
$$
\ll_\e T^{1/2+\e}K_0^{3/2},
$$
which is worse that the trivial estimation of $S(K)$, since (1.2) holds.
Likewise the use of the range of integration (5.21) gives also a poor bound.

\medskip
We shall conclude with a discussion on the shape of the functions
$\f_\ell(K,T;m,f)$, which appear in (2.1). We note that (see (3.9)) we have
$$
q_N(r) = 1 + \sum_{\ell=1}^L b_\ell r^{-2\ell} + O_{N,L}(r^{-2L-2})\leqno(5.25)
$$
with effectively computable constants $b_\ell$, where (as before) $L$ is taken
so large that the error term makes, in the appropriate expressions, a negligible
contribution. Each factor $r^{-2\ell}$ in (5.25) becomes, after change of variable
in the integral in (5.19),
$$
(K+Gu)^{-2\ell} = K^{-2\ell}\left\{1 + \sum_{j=1}^L d_\ell (Gu/K)^j +
O_{\ell,L}\left((Gu/K)^{L+1}\right)\right\},
$$
which is then evaluated by (5.3), furnishing a sum
containing powers of $G$ and $K$.

 In what concerns the factors $K^{1-\ell}(K-r)^\ell$ in (3.10),
note that$(K-r)^\ell$ introduces the factor $(Gu)^\ell$ in (5.19), and then
the corresponding integral is again evaluated by (5.3), producing eventually a suitable
power of $G$. The factor $K^{1-\ell}$, after the saddle point method is applied,
in view of (5.22) leads to
$$
K_1^{1-\ell} = (4T)^{1-\ell}(\sqrt{x}+\sqrt{1+x})^{2\ell-2}\qquad(x = m/f),
$$
and we have the power expansion (5.23). When this is all put together, we get
terms of the type $\f_\ell(K,T;m,f)$, which are power functions in each of
the variable, all of which are certainly $o(1)$ (as $K\to\infty$ and (2.1)
holds).

\bigskip
\head 6. Completion of proof of Theorem 1
\endhead
To complete the proof of Theorem 1 we shall show that
$$
{\Cal H}_4(f;h) \= \pi^{-3}\sum_{m=1}^{f-1}m^{-{1\over2}}d(m)d(f-m)
\Psi^-\left({m\over f};h\right)\leqno(6.1)
$$
makes a negligible contribution to (4.3). We use the representation
(this is [23, eq. (3.3.45)]), valid for $x = m/f < 1$ and $-1 < \b < -\hf$,
$$
\Psi^-(x;h) = \int\limits_0^\infty\Bigl\{\int\limits_{(\b)}
x^s(y(1+y))^{s-1}{\G^2(\hf-s)\d s\over\G(1-2s)\cos(\pi s)}\Bigr\}
\int\limits_{-\infty}^\infty rh(r)\left({y\over1+y}\right)^{ir}\d r\d y,
\leqno(6.2)
$$
where the triple integral converges absolutely. The function (6.2) can
be compared to the representation (5.11) for $\Psi^+(x;h)$: the function
$\Psi^-(x;h)$ is easier to deal with because of the factor $\cos(\pi s)$
in the denominator, and summation over $m$ in (6.1)
is finite. On the other hand, it has the drawback that the integral over $y$
is not finite, and there is an additional integration over $s$. As before,
it will suffice to consider the contribution of $|r\pm K|\le G\log K$.
Namely if $|r\pm K|\ge G\log K$ we interchange the order of integration,
and in the $y$ integral we integrate by parts the subintegral over
$(0,\,1]$ to obtain that the contribution is $\ll x^\b\exp(-\hf\log^2K)$.
For  $|r - K| \le G\log K$ (the case of the `+' sign is analogous) we make
the change of variable $r = K + Gu$ to obtain that the dominant contribution
of the $r$-integral will be
$$
GK{\roman e}^{iK\log{y\over1+y}}{\roman e}^{iK\log{4T\over K}}
\int\limits_{-\log K}^{\log K}\exp\Bigl(-u^2\pm iGu
\log{4T\over K}
+ iGu\log{y\over1+y}\Bigr)\d u.\leqno(6.3)
$$
Using (5.2) it follows that (6.3) becomes, up to a negligible error,
a multiple of
$$
\eqalign{&
GK\exp\left(iK\log\left({y\over1+y}\cdot{4T\over K}\right)\right)
\exp\left(-{\txt{1\over4}}G^2\left(\log\Bigl({y\over1+y}
\cdot{4T\over K}\Bigr)\right)^2
\right)\cr&
+ GK\exp\left(iK\log\left({y\over1+y}\cdot{K\over 4T}\right)\right)
\exp\left(-{\txt{1\over4}}G^2\left(\log
\Bigl({y\over1+y}\cdot{K\over 4T}\Bigr)\right)^2
\right).\cr&}\leqno(6.4)
$$
Since
$$
{\left(\log\Bigl({y\over1+y}\cdot{K\over 4T}\Bigr)\right)}^2
\ge \log^2\left({4T\over K}\right)\qquad(y > 0),
$$
this means that the contribution of the second exponential factor above
will be negligible, and the same holds for the first exponential factor,
if $y\ge 1$. In view of Stirling's formula (see (5.17)) and
$$
|\cos(x + iy)| \= \sqrt{\cos^2x + \sinh^2y}\qquad(x \in \RR,\, y \in \RR),
$$
it follows that the contribution of $|\I s| = |t| > \log^2K$ in (6.1)
will be negligibly small. If $0 \le y \le 1$ and
$$
{y\over1+y}\cdot{4T\over K} \le 1 - {\log T\over G}\leqno(6.5)
$$
or
$$
{y\over1+y}\cdot{4T\over K} \ge 1 + {\log T\over G},\leqno(6.6)
$$
the total contribution is again negligible. If (6.5) and (6.6) both fail,
then $y$ lies in an interval of length $\approx (K\log T)/(TG)$. But then  we
may integrate by parts the factor $y^{ir}$ in the integral, each time
increasing the exponent of $y$ by unity. If this is done sufficiently
many times, then trivial estimation shows that the total contribution
of (6.1) is negligibly small, and Theorem 1 is proved, if we take
$G = K_0^{1/2-\e}$ in (3.11) and (5.9) and replace $K_0$ by $K$.

\bigskip
\head 7. The proof of Theorem 2
\endhead
The proof of the first bound in (2.13) is straightforward. Namely
Motohashi derived the transformation formula for (4.1) by writing
$\H = H_j^2(\hf)\cdot H_j(\hf)$, and then by expressing $H_j(\hf)$
as a partial sum of $t_j(f)f^{-1/2}$ (see [23, Lemma 3.9]
or (7.7)) to which
the transformation formula for the bilinear sum of Hecke series is applied.
Therefore our problem reduces essentially to the evaluation and estimation
of Theorem 1 in the case $f = 1$. We obtain
$$
\eqalign{&
\sum_{K<\k_j\le K'<2K} \a_j H_j^2(\hf)\cos\left(\k_j\log\left(
{4{\roman e}T\over \k_j}\right)\right)\cr&
\asymp T\sum_{C_1TK^{-1}\le m
\le C_2TK^{-1}}
m^{-{3\over2}}d(m)d(m+1){\roman e}^{i{T\over m}} +
O_\e(K^{{3/2}+\e})\cr&
\ll_\e T^{1/2+\e}K^{1/2} +  K^{{3/2}+\e}
\ll_\e T^{1/2+\e}K^{1/2},\cr}\leqno(7.1)
$$
since (1.2) holds. We remark, similarly as in the discussion concerning
Theorem 1, that the sum over $m$ in (7.1) could be treated by the techniques
of [12]--[13] involving the binary additive divisor problem, but it seems
that the result that would be obtained in this fashion does not improve
the above (trivial) bound.

\medskip
For the proof of the second bound in (2.13) we proceed analogously to the
proof of
$$\sum_{\k_j\le T}\a_j H_j(\hf) = \left({T\over\pi}\right)^2
- BT\log T + O(T(\log T)^{1/2})\qquad(B>0), \leqno(7.2)
$$
given by M. Jutila and the author in [10].  The proof of (7.2) rested
on the use of (see e.g., [23] for a proof)

\medskip
{\bf Lemma 1.} (The first Bruggeman-Kuznetsov trace formula).
{\it Let $f(r)$ be an even, regular function for $|\I r| \le \hf$
such that $f(r) \ll (1+|r|)^{-2-\delta}$ for some $\delta>0$. Then}
$$\eqalign{&
\sum_{j=1}^\infty \a_j t_j(m)t_j(n)f(\k_j)
+ {1\over\pi}\int_{-\infty}^\infty {\s_{2ir}(m)\s_{2ir}(n)\over
(mn)^{ir}|\z(1+2ir)|^2}f(r)\d r\cr&
= {1\over\pi^2}\delta_{m,n}\int_{-\infty}^\infty r\tanh(\pi r)f(r)\d r
+ \sum_{\ell=1}^\infty{1\over\ell}
S(m,n;\ell)f_+\left({4\pi\sqrt{mn}\over\ell}\right),\cr}\leqno(7.3)
$$
{\it where $\delta_{m,n} = 1$ if $m = n$ and zero otherwise} ($m,n >0$),
$\s_a(d) = \sum_{d\mid n}d^a$, $S(m,n;\ell)$ {\it is the Kloosterman sum and}
$$
f_+(x) = {2i\over\pi}\int_{-\infty}^\infty {r\over\cosh(\pi r)}J_{2ir}(x)
f(r)\d r.\leqno(7.4)
$$

\medskip
In this formula one takes $n=1$ and $f(r) \equiv h_\ell(r;T,K,G)$, as given
by (3.10), and follows the scheme of proof of Theorem 1. This consists
of evaluating
$$\eqalign{&
{1\over\sqrt{\pi}G}\int\limits_{K_0}^{K_0'}
\sum_{|\k_j-K|\le G\log K_0}
\a_jH_j(\hf){\roman e}^{iK}\exp\left(i\k_j\log{4T\over K}\right)
{\roman e}^{-(\k_j-K)^2G^{-2}} \d K\cr&
= {1\over\sqrt{\pi}G}\int\limits_{K_0}^{K_0'}
\sum\nolimits_0(K;T,G){\roman e}^{iK}\d K
+ O(1),\cr}\leqno(7.5)
$$
where $G$ satisfies (3.2) and
$$
\sum\nolimits_0(K;T,G)
:= \sum_{j=1}^\infty\a_j H_j(\hf) h(\k_j;T,K,G).\leqno(7.6)
$$
To obtain the expression for (7.6) one multiplies  (7.3) by $m^{-1/2}$,
since (see [10] for proof) we have

\medskip
{\bf Lemma 2.} {\it Let $\k_j = (1+o(1))K,\, r = (1+ o(1))K$ } ($r\in \RR$)
{\it as
$K\to\infty,  Y = (1+\delta){K^2\over4\pi^2}$, with $\delta>0$
a given constant. Then, for any fixed positive constant $A>0$, there
exists a constant $C = C(A,\delta) > 0$ such that, for $h = C\log K$,
we have}
$$
H_j(\hf) = \sum_{m\le (1+\delta)Y}t_j(m)m^{-1/2}{\roman e}^{-(m/Y)^h}
+ O(K^{-A}),\leqno(7.7)
$$
{\it and}
$$
\z(\hf+ir)\z(\hf-ir) = \sum_{m\le (1+\delta)Y}\s_{2ir}(m)m^{-{1\over2}-ir}
{\roman e}^{-(m/Y)^h}+ O(K^{-A}).\leqno(7.8)
$$
In the proof of (7.2) the main term came from the integral
$$
\int_{-\infty}^\infty r\tanh(\pi r)f(r)\d r\leqno(7.9)
$$
in (7.3). However, now in the function $f(r)$ we shall have
the additional oscillating factor $(4T/K)^{\pm ir}$. Because of this,
 when we make the change of variable $r = \pm K + Gu$,
we shall eventually wind up with exponential factors of the form
$$
\exp\left\{-{\txt{1\over4}}G^2\left(\log{4T\over K}\right)^2\right\},
$$
which make a negligible contribution. The total contribution of the
continuous spectrum (the integral on the left-hand side of (7.3))
is easily seen to be $\ll_\e K_0^{1+\e}$. The only delicate part
is the Kloosterman-sum contribution, coming from the right-hand side
of (7.3). However, this presents no major problem, since the
estimation is  analogous to the one made in [10] for the
proof of (7.2). We shift
the line of integration in the integral defining $f_+$ to
$\I r = -1$ and use the power series representation
$$
J_{2+ix}(z) = \sum_{k=0}^\infty {(-1)^k(z/2)^{2+ix+2k}
\over\G(k+1)\G(k+ 2+ix+1)} \quad (z = 4\pi\sqrt{m}/\ell \ll K^{1-B}),
$$
which shows that the contribution of $\ell > K^{B}$ is $\ll K^{-A}$
for any fixed $A>0$, provided that $B = B(A)$ is sufficiently large.
The only difference from [10] is that, in making the shift, the factor
$(4T/K)^{ir}$ will make now a contribution of $O(T/K)$, which is harmless
if $B$ is sufficiently large.
In the remaining sum, we substitute (see e.g., [18, p. 139])
$$
J_{2ir}(x) - J_{-2ir}(x) = {2i\over\pi}\sinh(\pi r)
\int_{-\infty}^\infty \cos(x\cosh u)\cos(2ru)\d u.
$$
Integration by parts shows that, for $x>0$ and $r\ge0$,
$$
\eqalign{
J_{2ir}(x) - J_{-2ir}(x)
& = {2i\over\pi}\sinh(\pi r)
\int_{-\log^2K}^{\log^2K} \cos(x\cosh u)\cos(2ru)\d u\cr&
+ O\left(x^{-1}(r+1)\exp(\pi r - \hf\log^2K)\right).\cr}\leqno(7.10)
$$
The error term in (7.10) clearly makes a negligible contribution.
 The main term in (7.10) will contribute to $f_+$
$$
-{4\over\pi^2}\int_{-\log^2K}^{\log^2K} \cos(x\cosh u)\int_0^\infty
rf(r,K)\tanh(\pi r)\cos(2ru)\d r\d u,\leqno(7.11)
$$
where
$$
x \;=\; 4\pi{\sqrt{m}\over\ell} \le 2(1+\delta)K.\leqno(7.12)
$$
In the inner integral we use (5.12) and make
the change of variable $r = K + Gv$. In the ensuing $v$-integral the
non-negligible contribution will be from the range $|v| \le \log K$.
Since $f(r)$ contains the factor $(4T/K)^{ir}$, it follows by (5.2)
and (5.3) that the contribution of $f_+$ is
$$
\asymp
\R\left\{GK\int_{-\log^2K}^{\log^2K} \cos(x\cosh u)
\exp\left(-{G^2\over 4}\bigl(\log{4T\over K} \pm 2u\bigr)^2
\pm 2iKu)\right)\d u\right\}.\leqno(7.13)
$$
The relevant exponential factor  will be of the form
$$
\exp(ig(u)),\;g(u) = x\cosh u \pm 2Ku,\; g'(u) = x\sinh u \pm 2K.
$$
The saddle point $u_1$ is (here the solution of $g'(u_1) = 0$
with the plus sign is treated, since  the other case is similar)
$$
u_1 = \log\left({2K\over x} + \sqrt{{4K^2\over x^2}+1}\,\right),
$$
and we have
$$
g''(u_1) = x\cosh(u_1) \gg K.
$$
Since $K/x \gg 1$ in view of (7.12), it follows by the saddle
point method that the main contribution to (7.11) is
$$
\asymp \int_{K_0}^{K_0'}{\roman e}^{\pm iK + iH(K)}
K^{1/2}\exp\left(-{G^2\over4}{\left(
\log {4T/K\over 2K/x + \sqrt{(2K/x)^2+1}}\right)}^2\right)\d K,
\leqno(7.14)
$$
plus an error term which does not exceed $O(T^{1/2+\e}K^{1/4})$,
where
$$
H(K) := g(u_1),\quad |H'(K)| = \log \left({2K\over x}
 + \sqrt{{4K^2\over x^2}+1}\,\right) + O(1),
$$
and the contribution is negligible unless
$$
{C_1T\over K_0^2}\sqrt{m} \le \ell \le {C_2T\over K_0^2}\sqrt{m}
\qquad(0 < C_1 < C_2).\leqno(7.15)
$$
Thus by the first derivative test the integral in (7.14) is
$\ll K_0^{1/2}\log K_0$.
If we use Weil's classical bound $|S(m,n;\ell)| \le
(m,n,\ell)^{1/2}d(\ell)\ell^{1/2}$, then we see that the total contribution
of the Kloosterman sum term in (7.3) is
$$
\eqalign{&
\ll_\e K_0^{1/2+\e}\sum_{m\ll K}m^{-1/2}
\sum_{\ell \approx {T\over K_0^2}\sqrt{m}}
{1\over \ell}|S(m,1;\ell)| \cr&
\ll_\e K_0^{1/2+\e}\sum_{m\ll K_0}m^{-1/2}
\sum_{\ell \approx {T\over K_0^2}\sqrt{m}}d(\ell)\ell^{-1/2}\cr&
\ll_\e K_0^{\e-1/2}T^{1/2}\sum_{m\ll K_0}m^{-1/4}\cr&
\ll_\e T^{1/2+\e}K_0^{1/4}.\cr}
$$
We take  $G = K_0^{\e}$, note that $K_0 \ll T^{1/2}K_0^{1/4}$ in view
of (1.2) and finally replace $K_0$ by $K$. Then the second bound in
(2.13) follows and  the proof of Theorem 2 is complete.

\bigskip
\head 8. Proof of Theorem 3
\endhead
Suppose that the hypotheses of Theorem 3 hold. We start from
$$
\int_T^{2T}{\bigl(S_m(K;K',t)\bigr)}^2\d t
\le \int_{T/2}^{5T/2}\f(t){\bigl(S_m(K;K',t)\bigr)}^2\d t,\leqno(8.1)
$$
where $\f(t)$ is a non-negative, smooth function supported in $[T/2,\,5T/2]\,$
such that $\f(t) = 1$ for $T\le t \le 2T$. We assume that $m = 3$, as this is
the most interesting case. The proof of the cases $m = 1,2$ is analogous, only
instead of (1.4)--(1.5) we shall need the corresponding
bounds with $H_j^2(\hf)$ (see [23, eq. (3.4.4)]) or $H_j(\hf)$ (see [10]).
If the cosine is written as a sum of exponentials, then for $m=3$ the
right-hand side of (8.1) becomes, after integration by parts,
$$
\eqalign{&
\ll \int\limits_{T/2}^{5T/2}\f(t)\sum_{K<\k_j,\k_\ell\le K'}
\a_j\a_\ell H_j^3(\hf)H_\ell^3(\hf){\roman e}^{i(\k_\ell\log\k_\ell-\k_j\log\k_j)}
\left(4{\roman e}t\right)^{i\k_j-i\k_\ell}\d t\cr&
= -\sum_{K<\k_j,\k_\ell\le K'}
\a_j\a_\ell H_j^3(\hf)H_\ell^3(\hf)
{\roman e}^{i(\k_\ell\log\k_\ell-\k_j\log\k_j)}\cr&\times
\int\limits_{T/2}^{5T/2}{\f'(t)\over i\k_j-i\k_\ell+1}
\left(4{\roman e}\right)^{i\k_j-i\k_\ell}t^{i\k_j-i\k_\ell+1}\d t.
\cr}\leqno(8.2)
$$
In (8.2) we may continue to integrate by parts, noting that
$$
\f^{(r)}(T/2) = \f^{(r)}(5T/2) = 0,\quad
\f^{(r)}(t) \;\ll_r\; T^{-r}\qquad(r = 0,1,2,\ldots\;).\leqno(8.3)
$$
Therefore taking $r = r(A,\e)$ sufficiently large, it follows from (8.3) that the
contribution of $\k_j,\k_\ell$ such that $|\k_j-\k_\ell| > T^\e$ is $\ll T^{-A}$
for any given, large $A>0$. The contribution of the remaining pairs
$\k_j,\k_\ell$ is estimated trivially by the use of (1.3)--(1.5) as
$$
\eqalign{&
\ll \int\limits_{T/2}^{5T/2}\f(t)\sum_{K<\k_j\le K'}\a_j H_j^3(\hf)
\sum_{|\k_j-\k_\ell| \le T^\e}\a_\ell H_\ell^3(\hf)\d t\cr&
\ll_\e T^\e K\int\limits_{T/2}^{5T/2}\f(t)
\sum_{K<\k_j\le K'}\a_j H_j^3(\hf)\d t \ll_\e T^{1+\e}K^3,
\cr}
$$
and this is asserted by (2.15). If the conjectural (1.7)--(1.8) holds, then
obviously (2.15) can be improved (for $m=3$) to
$$
\int_T^{2T}{\bigl(S(K;K',t)\bigr)}^2\d t \,\ll_\e\, T^{1+\e}K^{5/2}.
$$
Also by direct integration we have
$$
\int_T^{2T}S(K;K',t)\d t \,\ll_\e\, T^{1+\e}K,\leqno(8.4)
$$
while the integral in (8.4) is $\ll_\e T^{1+\e}K^{1/2}$ if (1.7)--(1.8) holds.
\medskip
Finally we sketch the proof of (2.16) of the Corollary. We start from
$$
\int_T^{2T}{(E_2(2t)-E_2(t))}^2\d t
\le \int_{T/2}^{5T/2}\f(t){(E_2(2t)-E_2(t))}^2\d t,\leqno(8.5)
$$
where $\f(t)$ is as in (8.1).
Then we use (2.7)--(2.8), truncating the series in (2.8) at $T\D^{-1}\log T$
with a negligible error. After this, we remove the monotonic coefficients
$\k_j^{-3/2}$ and $\exp\Bigl(-{\txt{1\over4}}
\bigl({\D\k_j\over T}\bigr)^2\Bigr)$ by partial summation. Then we
obtain  the sum $S_m(K;K',t)$ with $m=3$ and $t$ replaced by
$2t+\D\log T$ or $t-\D\log T$, which does not cause any trouble. Hence the
integral on the left-hand side of (8.5) is
essentially majorized by $\ll_\e T^\e$ integrals of the type
$$
T\int_{T/2}^{5T/2}\f(t){\bigl(K^{-3/2}S_m(K;K',t)\bigr)}^2\d t
\ll_\e T^{2+\e},
$$
and (2.16) follows on replacing $t$ by $t2^{-j}$ in the integrand in
(8.5), and summing up the corresponding bounds over $j = 1,2,\ldots\,$.

\medskip
It may be remarked that the method of proof of Theorem 3 gives also,
for $1 \ll K < K' \le 2K \ll T$,
$$
\int_T^{2T}{\bigl(S_m(K;K',t)\bigr)}^4\d t
\;\ll_\e\; T^{1+\e}K^7\qquad(m = 1,2,3),
$$
which means that, in the mean fourth sense, the sum
$S_m(K;K',t)$ is $\ll_\e K^{7/4+\e}$.

\bigskip
\Refs
\bigskip

\item{[1]}  D.R. Heath-Brown, The fourth moment of the Riemann
zeta-function,
{\it Proc. London Math. Soc.}  {\bf(3)38}(1979), 385-422.

\item {[2]} A. Ivi\'c,  The Riemann zeta-function, {\it John Wiley and
Sons}, New York, 1985 (2nd ed. {\it Dover}, 2003).

\item {[3]} A. Ivi\'c,  Mean values of the Riemann zeta-function, LN's
{\bf 82}, {\it Tata Institute of Fundamental Research}, Bombay,  1991
(distr. by Springer Verlag, Berlin etc.).

\item{[4]}  A. Ivi\'c,  On the fourth moment of the Riemann
zeta-function, {\it Publs. Inst. Math. (Belgrade)} {\bf 57(71)}
(1995), 101-110.

\item{[5]} A. Ivi\'c,  The Mellin transform and the Riemann
zeta-function,  {\it Proceedings of the Conference on Elementary and
Analytic Number Theory  (Vienna, July 18-20, 1996)},  Universit\"at
Wien \& Universit\"at f\"ur Bodenkultur, Eds. W.G. Nowak and J.
Schoi{\ss}engeier, Vienna 1996, 112-127.

\item{[6]} A. Ivi\'c, On the error term for the fourth moment of the
Riemann zeta-function, {\it J. London Math. Soc.},
{\bf60}(2)(1999), 21-32.

\item {[7]} A. Ivi\'c, On sums of Hecke series in short intervals,
{\it Journal de Th\'eorie des Nombres de Bordeaux} {\bf13}(2001), 453-468.

\item {[8]} A. Ivi\'c, On some conjectures and results for the Riemann
zeta-function and Hecke series, Acta Arith. {\bf99}(2001), 115-145.

\item{[9]} A. Ivi\'c,  On the moments of Hecke series at central points,
{\it Funct. Approximatio} {\bf30}(2002), 49-82.

\item{[10]} A. Ivi\'c and M. Jutila,
On the moments of Hecke series at central points II,
{\it Funct. Approximatio} {\bf31}(2003), 7-22.

\item{ [11]} A. Ivi\'c and Y. Motohashi,  The mean square of the error
term for the fourth moment of the zeta-function,  {\it Proc. London
Math. Soc.} (3){\bf 66}(1994), 309-329.

\item {[12]} A. Ivi\'c and Y. Motohashi,  The fourth moment of the
Riemann zeta-function, {\it J. Number Theory} {\bf 51}(1995), 16-45.

\item {[13]}   A. Ivi\'c and Y. Motohashi, On some estimates involving
the binary additive divisor problem, {\it Quart. J. Math. (Oxford)}
(2){\bf46}(1995), 471-483.

\item {[14]} M. Jutila, On exponential sums involving the divisor
function, {\it Journ.  reine angew. Math.} {\bf355}(1985), 173-190.

\item {[15]} S. Katok and P. Sarnak, Heegner points, cycles and Maass
forms, {\it Israel J. Math.} {\bf84}(1993), 193-227.

\item {[16]} N.V. Kuznetsov, On the eigenfunctions of an
integral equation (in Russian),
{\it Zapiski Nauchnykh Seminarov LOMI} {\bf17}(1970), 66-149.

\item {[17]} N.V. Kuznetsov, Convolution of the Fourier
coefficients of the Eisenstein--Maass series (in Russian),
{\it Zapiski Nauchnykh Seminarov LOMI} {\bf129}(1983), 43-84.

\item {[18]} N.N. Lebedev, Special functions and their applications,
{\it Dover Publications, Inc.}, New York, 1972.

\item {[19]} Y. Motohashi, Spectral mean values of Maass wave forms,
{\it J. Number Theory} {\bf 42}(1992), 258-284.

\item{ [20]} Y. Motohashi,   An explicit formula for the fourth power
mean of the Riemann zeta-function, {\it Acta Math. }{\bf 170}(1993),
181-220.

\item{ [21]} Y. Motohashi, The binary additive divisor problem,
{\it Annales Scien. \'Ecole Norm. Sup.}, $4^e$ s\'erie,
{\bf27}(1994), 529-572.

\item {[22]} Y. Motohashi,  A relation  between the Riemann
zeta-function and the hyperbolic Laplacian, {\it Annali Scuola Norm.
Sup. Pisa, Cl. Sci. IV ser.} {\bf 22}(1995), 299-313.

\item{[23]} Y. Motohashi, Spectral  theory of the Riemann zeta-function,
{\it Cambridge University Press}, Cambridge, 1997.

\item {[24]} N.I. Zavorotnyi, On the fourth moment of the Riemann
zeta-function (in Russian), Automorphic functions and number
theory I, {\it Collected Scientific Works}, Vladivostok, 1989, 69-125.

\vskip5cm
\endRefs

\bigskip
\sevenrm
Aleksandar Ivi\'c

Katedra Matematike RGF-a

Universitet u Beogradu, \DJ u\v sina 7

 11000 Beograd, Serbia and Montenegro

{\sevenbf ivic\@rgf.bg.ac.yu}

\vfill


\bye